\documentclass[12pt]{amsart}

\date{\today}

\usepackage[margin=1.15in]{geometry}

\usepackage{amsmath,amscd,amssymb,amsfonts,latexsym,wasysym, mathrsfs, enumitem, mathtools,hhline,color}
\usepackage[all, cmtip]{xy}

\usepackage{url}

\definecolor{hot}{RGB}{65,105,225}

\usepackage[pagebackref=true,colorlinks=true, linkcolor=hot ,  citecolor=hot, urlcolor=hot]{hyperref}

\usepackage{graphicx,enumerate}

\newcommand{\CN}{\mathbb{C}^{n}}
\newcommand{\C}{\mathbb{C}}
\newcommand{\N}{\mathbb{N}}
\newcommand{\R}{\mathbb{R}}

\theoremstyle{plain}
\newtheorem{theorem}{Theorem}[section]
\newtheorem{prop}[theorem]{Proposition}

\newtheorem{lm}[theorem]{Lemma}

\newtheorem{cor}[theorem]{Corollary}

\newtheorem{thrm}[theorem]{Theorem}

\theoremstyle{definition}

\newtheorem{defn}[theorem]{Definition}
\newtheorem{que}[theorem]{Question}
\newtheorem{rmk}[theorem]{Remark}

\newtheorem{ex}[theorem]{Example}
\newtheorem*{ex*}{Example}

\def\be{\begin{equation}}
\def\ee{\end{equation}}

\def\bt{\begin{thrm}}
\def\et{\end{thrm}}

\def\bc{\begin{cor}}
\def\ec{\end{cor}}

\def\br{\begin{rmk}}
\def\er{\end{rmk}}

\def\bp{\begin{prop}}
\def\ep{\end{prop}}

\def\bl{\begin{lm}}
\def\el{\end{lm}}

\def\bex{\begin{ex}}
\def\eex{\end{ex}}

\def\bd{\begin{defn}}
\def\ed{\end{defn}}

\DeclareMathOperator{\kn}{Ker}

\title{Milnor and Tjurina numbers for a hypersurface germ with isolated singularity}

\author{Yongqiang Liu}

\address{KU Leuven, Department of Mathematics,
Celestijnenlaan 200B, 3001 Leuven, Belgium} 

\email{liuyq1117@gmail.com}

\date{\today}



\begin{document}

\maketitle

\begin{abstract}  Assume that $f:(\CN,0) \to (\C,0)$ is an analytic function germ at the origin with only isolated singularity. Let $\mu$  and $\tau$ be the corresponding Milnor and Tjurina numbers. We show that $\dfrac{\mu}{\tau} \leq n$.  As an application, we give a lower bound for the Tjurina number in terms of $n$ and the multiplicity of $f$ at the origin. 
 \end{abstract}
 
\section{Main result}

Assume that $f:(\CN,0) \to (\C,0)$ is an analytic function germ at the origin with only isolated singularity. Set $X=f^{-1}(0)$.
 Let $S=\C\lbrace x_{1}, \ldots, x_{n} \rbrace$ denote the formal power series ring. Set $J_{f}=(\partial f/\partial x_{1}, \ldots, \partial f/ \partial x_{n})$ as the Jacobian ideal. Then the Milnor and   Tjurina algebras are defined as  \begin{center}
  $M_{f}=S/J_{f}$, and  $T_{f}= S /(J_{f},f)$.
 \end{center}
  Since $X$ has isolated singularities, $M_{f}$ and $T_{f}$ are finite dimensional $\C$-vector spaces. The corresponding dimension $\mu$  and $\tau$ are called the Milnor and Tjurina numbers, respectively. It is clear that $\mu \geq \tau$.

Consider the following long exact sequence of $\C$-algebras: 
\be
0 \to \kn(f) \to M_{f} \overset{ f}{\to} M_{f} \to T_{f} \to 0
\ee
where the middle map is  multiplication by $f$, and  $\kn(f)$ is the kernel of this map. Then $\dim_{\C} \kn(f) = \tau$. 

Recall a well known result given by J. Briancon and H. Skoda in \cite{BS}:  $$f^{n} \in J_{f},$$ which shows that $f^{n}=0$ in $M_{f}$, i.e., $(f^{n-1}) \subset \kn(f)$. Here $(f^{n-1})$ is the ideal in $M_{f}$ generated by $f^{n-1}$. The following theorem is a direct application of this result.  
  
\bt \label{0.1} Assume that $f:(\CN,0) \to (\C,0)$ is an analytic function germ at the origin with only isolated singularity. Then $$ \dfrac{\mu}{\tau} \leq n.$$ Moreover, $\dfrac{\mu}{\tau}= n$, if and only if, $\kn(f)=(f^{n-1})$.  
\et 
\begin{proof}
 Since $f^{n}=0$ in $M_{f}$, we have the following finite decreasing filtration: \begin{center}
$ M_{f} \supset (f) \supset (f^{2})  \supset \cdots \supset (f^{n-1}) \supset (f^{n})=0$
\end{center}  where $(f^{i})$ is the ideal in $M_{f}$ generated by $f^{i}$.

Consider the following long exact sequence:
\be 
0 \to \kn(f)\cap (f^{i}) \to (f^{i}) \overset{ f}{\to} (f^{i}) \to (f^{i})/(f^{i+1}) \to 0
\ee  
 where the middle map is multiplication by $f$. Then,
 $$ \dim_{\C} \lbrace(f^{i})/(f^{i+1})\rbrace=\dim_{\C}\lbrace \kn(f)\cap (f^{i})\rbrace\leq \dim_{\C} \kn(f)= \tau . $$  
  Therefore, $$\mu= \dim_{\C} M_{f} =\dim_{\C} T_{f}+ \sum_{i=1}^{n-1} \dim_{\C}\lbrace (f^{i})/(f^{i+1})\rbrace 
  \leq n \cdot \tau .$$
  
 $\dfrac{\mu}{\tau}= n$, if and only if, for any $1\leq i \leq n-1$, $\kn(f)\cap (f^{i})= \kn(f)$, i.e., $\kn(f) \subset (f^{i})$. On the other hand, $(f^{n-1}) \subset \kn(f)$. Hence, $\kn(f) = (f^{n-1})$.   
\end{proof}

 K. Saito showed (\cite{S}) that $\dfrac{\mu}{\tau}=1$ holds, if and only if, $f$ is weighted homogeneous, i.e., analytically equivalent to such a polynomial. It leads to the following natural question.
\begin{que} Is this upper bound of $\dfrac{\mu}{\tau}$ optimal?  When can the optimal upper bound be obtained?
\end{que}
\br  Recently,  A. Dimca and G.-M. Greuel showed (\cite[Theorem 1.1]{DG}) that the upper bound $\dfrac{\mu}{\tau}\leq 2$ can never be achieved for the isolated plane curve singularity case unless $f$ is smooth at the origin. Moreover, they gave (\cite[Example 4.1]{DG}) a  sequence of isolated plane curve singularity with the ratio $\dfrac{\mu}{\tau}$ strictly increasing towards  $4/3$. In particular, the singularities can be chosen to be all either irreducible, or consisting of smooth branches with distinct tangents. Based on these computations, they asked (\cite[Question 4.2]{DG}) whether 
$$ \dfrac{\mu}{\tau} < 4/3$$
for any isolated plane curve singularity?
\er
\bex It is clear that $\dfrac{\mu}{\tau}>n-1$ implies that $f^{n-1} \notin J_{f}$. 

  Consider the function germ: $$ f=(x_{1}\cdots x_{n})^{2}+ x_{1}^{2n+2} +\cdots +x_{n}^{2n+2},$$ which defines an isolated singularity at the origin.  B. Malgrange showed (\cite{M}) that the monodromy on the $(n-1)$-th cohomology of the Milnor fibre has a Jordan block with size $n$. Coupled with the theorem by J. Scherk (\cite[Theorem]{Sch}), it gives us that $f^{n-1} \notin J_{f}$. It can be checked with the software SINGULAR that $\dfrac{\mu}{\tau}<1.5$ for $n\leq 7$, which is far away from our upper bound $n$. 
\eex

\section{Applications}
Theorem \ref{0.1} implies a well known result in complex singularity theory, which states that the Milnor number of an analytic function germ  is finite (or non-zero) if and only if the Tjurina number is so (see \cite[Lemma 2.3, Lemma 2.44]{GLM}). 

\subsection{A lower bound for the Tjurina number} First we recall a well-known lower bound for $\mu$ in terms of $n$ and the multiplicity $m$ of $f $ at the origin. The following description can be found in \cite{Huh}.

The sectional Milnor numbers associated to the germ $X$ are introduced by Teissier \cite{T1}. The $i$-th sectional Milnor number of the germ $X$, denoted $\mu^{i}$, is the Milnor number of the intersection of $X$ with a general $i$-dimensional plane passing through the origin (it does not depend on the choice of the generic planes). Then $\mu=\mu^{n}$.  The Minkowski inequality for mixed multiplicities says that the sectional Milnor numbers
always form a log-convex sequence \cite{T2}. In other words, we have
$$\dfrac{\mu^{n}}{\mu^{n-1}}  \geq \dfrac{\mu^{n-1}}{\mu^{n-2}}\geq \cdots \geq \dfrac{\mu^{1}}{\mu^{0}} ,$$
where $\mu^{0}=1$ and $\mu^{1}=m-1$.  Then 
\be \label{2.1}
\mu \geq (m-1)^{n}.
\ee  Moreover,  the equality holds if and only if $f$ is a semi-homogeneous function (i.e., $f=f_{m}+g$, where $f_{m}$ is a homogeneous polynomial of degree $m$ defining an isolated singularity at the origin and $g$ consists of terms of degree at least $m+1$) after a biholomorphic change of coordinates.  For a detailed proof, see \cite[Proposition 3.1]{YZ}.

 The next corollary is a direct consequence  of Theorem \ref{0.1} and (\ref{2.1}). 
\bc Assume that $f:(\CN,0) \to (\C,0)$ is an analytic function germ at the origin with only isolated singularity. Then  $$\tau \geq \dfrac{(m-1)^{n}}{n} .$$
\ec
It is clear that, even for homogeneous polynomial case, this lower bound can never be obtained when $n>1$. In fact, in this case, $\tau =\mu =(m-1)^{n} >  \dfrac{(m-1)^{n}}{n}$.

\subsection{Another lower bound for the Tjurina number} Another lower bound for $\mu$ is given by A. G. Kushnirenko using the Newton number (\cite{K}). Let $\Gamma$ be the boundary of the Newton polyhedron of $f$, i.e., $\Gamma$ is a polyhedron of dimension $n-1$ in $\N^{n}$ (where $\N=\lbrace 0,1,2, \cdots \rbrace$) determined in the usual way by the non-zero coefficients in $f$. Then $f$ is said to be convenient if $\Gamma$ meets each of the coordinate axes of $\R^{n}$. 
Let $S$ be the union of all line segments in $\R^n$ joining the origin to points of $\Gamma$. For a convenient $f$, the Newton number $\nu(f)$ is defined as: $$ \nu=n!V_{n}-(n-1)!V_{n-1}+\cdots +(-1)^{n-1} 1! V_{1} +(-1)^{n}.　$$
where $V_n$ is the $n$-dimensional volume of $S$ and for $1\leq q \leq n-1$, $V_q$ is the sum of the $q$-dimensional volumes of the intersection of $S$ with the coordinate planes of dimension $q$.
 A. G. Kushnirenko showed that, if $f$ is convenient, then $$\mu \geq \nu,$$
 Moreover,  $\mu=\nu$ holds, if $f$ is non-degenerate. (For the definition of non-degenerate, see \cite[Definition 1.19]{K}.)  Again, this gives us a corresponding lower bound for the Tjurina number.
\bc Assume that $f:(\CN,0) \to (\C,0)$ is an analytic function germ at the origin with only isolated singularity, which is convenient. Then  $$\tau \geq \dfrac{\nu}{n} ,$$ where $\nu$ is the Newton number.
\ec
 \begin{que} Are the lower bounds of $\tau$ in Corollary 2.1 and 2.2 optimal? When can the optimal lower bounds be obtained?
\end{que}

For some special class of polynomials, the bound for the ratio $\dfrac{\mu}{\tau}$ can be improved. For example, A. Dimca showed that $f^{2} \in J_{f}$ for semi-weighted homogeneous polynomials (\cite[Example 3.5]{D}), hence $\dfrac{\mu}{\tau} \leq 2$ and $\tau \geq   \dfrac{(m-1)^{n}}{2}$ in this case.

 \bex   Choose $f=x^{m}+y^{m}+z^{m}+g$, where $g$ has degree at least $m+1$. Then $\mu=(m-1)^{3}$. It is shown in \cite[Example 4.7]{W}  that $\tau_{min}=(2m-3)(m+1)(m-1)/3$, when $g$ varies. 
 \eex
 
\textbf{Acknowledgments.} The author is grateful to Alexandru Dimca for useful discussions.  The author is partially supported by Nero Budur's research project G0B2115N from the Research Foundation of Flanders.


\end{document}